\documentclass{article}
\newtheorem{thm}{Theorem}[section]
\newtheorem{cor}[thm]{Corollary}
\newtheorem{lem}[thm]{Lemma}
\newtheorem{de}[thm]{Definition}
\newcounter{bean}
\newcounter{milk}

\begin{document}

\title{The Shape Parameter in the Shifted Surface Spline}                  
\author{ Lin-Tian Luh \\Department of Mathematics, Providence University\\ Shalu, Taichung\\ email:ltluh@pu.edu.tw \\ phone:(04)26328001 ext. 15126 \\ fax:(04)26324653   }        
\date{\today}          
\maketitle
{\bf Abstract.} The purpose of this article is to explore the optimal choice of shape parameter which is an important and longstanding problem in the theory of radial basis functions(RBF). We already handled it for multiquadric and Gaussian in \cite{Lu3,Lu4,Lu5,Lu6,Lu7}. Here we focus on shifted surface spline and present concrete criteria for the choice of shape parameter.\\
\\
{\bf key words}: radial basis function, shifted surface spline, shape parameter.\\
\\
{\bf AMS subject classification}: 41A05, 41A15, 41A25, 41A30, 41A63, 65D10
\section{Introduction}       
The so-called shifted surface spline is defined by
\begin{equation}
  h(x):=(-1)^{m}(|x|^{2}+c^{2})^{\frac{\lambda}{2}}log(|x|^{2}+c^{2})^{\frac{1}{2}},\ \lambda\in Z_{+},\ m=1+\frac{\lambda}{2},\ c>0,\ x\in R^{n},\ \lambda,\ n\ \mbox{even}
\end{equation}
, and
\begin{equation}
  h(x):=(-1)^{\lceil \lambda -\frac{n}{2}\rceil}(|x|^{2}+c^{2})^{\lambda-\frac{n}{2}},\ \lambda\in Z_{+}=\{ 1,2,3,\cdots \}, \ n\ \mbox{odd}\ \mbox{and}\ \lambda>\frac{n}{2}
\end{equation}
, where $|x|$ is the Euclidean norm of $x$, and $\lambda,\ c$ are constants. The number $\lceil\lambda-\frac{n}{2}\rceil$ denotes the smallest integer greater than or equal to $\lambda-\frac{n}{2}$.

This function was introduced by Dyn, Levin and Rippa in \cite{Dy1} for $x\in R^{2}$. Then it's extended to $R^{n}$ for $n\geq 1$. For further details we refer the reader to \cite{Bu,Du,Dy2,Po,Ri,Yo1,Yo2,Yo3,Yo4}.

Note that (2) is just multiquadric and will not be explored here. In this paper we only deal with (1). As is well known, $h(x)$ is conditionally positive definite of order $m=1+\frac{\lambda}{2}$ where $\lambda$ is of our choice. The constant $c$ is just the so-called shape parameter whose optimal choice is of our primary concern.

Since $h(x)$ is c.p.d.(conditionally positive definite) of order $m$, there exists a unique function $s(x)$ of the form
\begin{equation}
  s(x)=p(x)+\sum_{j=1}^{N}c_{j}h(x-x_{j})
\end{equation}
which interpolates any scattered set of data points $(x_{1},y_{1}),\cdots , (x_{N},y_{N})$ as long as $\{ x_{1},\cdots , x_{N}\}$ is polynomial-nondegenerate in $R^{n}$. Here $p(x)\in P_{m-1}$ is a polynomial of degree $\leq m-1$ and $c_{j}'s$ are chosen so that
$$\sum_{j=1}^{N}c_{j}q(x_{j})=0$$
for all polynomials $q$ in $P_{m-1}$.

We will show that whenever $x$ and the data points are fixed, the upper bound of $|f(x)-s(x)|$ will be made minimal by a good choice of $c$. Here $f(x)$ denotes the underlying function of the data points.

\subsection{Basic Ingredients}    
\begin{de}
For $n=1,2,\cdots ,$ the integers $\gamma_{n}$ are defined by $\gamma_{1}=2$ and $\gamma_{n}=2n(1+\gamma_{n-1})$ for $n>1$. 
\end{de}

The numbers $\gamma_{n}$ will appear in our criteria.

In this paper all approximated functions belong to a semi-Hilbert space ${\cal C}_{h,m}(R^{n})$, abbreviated as ${\cal C}_{h,m}$, where $m$ denotes the order of conditional positive definiteness of the function $h$. There is a semi-norm on ${\cal C}_{h,m}$, denoted by $\| \cdot \|_{h}$. The definition and characterization of this space can be found in \cite{Lu1,Lu2,MN1,MN2}. We will not repeat them. A subspace of ${\cal C}_{h,m}$, denoted by $B_{\sigma}$, is of our special concern and is defined as follows. 
\begin{de}
For any $\sigma>0$,
$$B_{\sigma}:=\{ f\in L^{2}(R^{n}):\ \hat{f}(\xi)=0\ for \ all\ |\xi|>\sigma \}$$
is the set of band-limited functions, where $\hat{f}$ denotes the Fourier transform of $f$.
\end{de}

We need some constants as follows.
\begin{de}
Let $n,\lambda$ and $m$ be as in (1). The constants $\rho$ and $\Delta_{0}$ are defined according to the following three cases.
\begin{list}
  {(\alph{bean})}{\usecounter{bean} \setlength{\rightmargin}{\leftmargin}}
  \item Suppose $n-\lambda >3$. Let $s=\lceil \frac{n-\lambda -3}{2}\rceil $. Then 
      $$\rho=1+\frac{s}{2m+3} \ and \ \Delta_{0}=\frac{(2m+2+s)(2m+1+s)\cdots (2m+3)}{\rho^{2m+2}}.$$     
  \item Suppose $n-\lambda \leq 1$. Let $s=-\lceil \frac{n-\lambda-3}{2}\rceil$. Then
   $$\rho=1\ and\ \Delta_{0}=\frac{1}{(2m+2)(2m+1)\cdots (2m-s+3)}.$$
  \item Suppose $1<n-\lambda\leq 3$. Then
 $$\rho =1\ and \ \Delta_{0}=1.$$  
\end{list}
\end{de}

Finally, the Fourier transform \cite{GS} of (1) should also be introduced. It's of the form
\begin{equation}
\hat{h}(\xi)=l(\lambda,n)|\xi|^{-\lambda-n}\tilde{{\cal K}}_{\frac{n+\lambda}{2}}(c|\xi|)
\end{equation}
where $l(\lambda,n)=(2\pi)^{-\frac{n}{2}}\cdot 2^{\frac{\lambda}{2}}\cdot (\frac{\lambda}{2}!)$ \cite{MN2} and $\tilde{{\cal K}}_{\nu}(t)=t^{\nu}{\cal K}_{\nu}(t)$, ${\cal K}_{\nu}(t)$ being the modified Bessel function of the second kind \cite{AS}.
\section{Fundamental Theory}
The cornerstone of our theory is the exponential-type error bound raised by the author in \cite{Lu8}. We cite it directly with a slight modification to make it easier to understand.

In the following theorem
$$d=d(E,X):=\sup_{y\in E}\inf_{x\in X}|y-x|$$
denotes fill distance.
\begin{thm}
Suppose $h$ is defined as in (1). Then, given any positive number $b_{0}$, there are positive constants $d_{0}$ and $\omega,\ 0<\omega <1$, which depend on $b_{0}$, for which the following is: If $f\in {\cal C}_{h,m}$ and $s$ is the $h$ spline defined in (3) that interpolates $f$ on a subset $X$ of $R^{n}$, then
\begin{equation}
|f(x)-s(x)|\leq \sqrt{l(\lambda,n)}\cdot (2\pi)^{\frac{1}{4}}\cdot \sqrt{n\alpha_{n}}\cdot \sqrt{\Delta_{0}}\cdot c^{\frac{\lambda}{2}}\cdot\omega^{\frac{1}{d}}\cdot \|f\|_{h}
\end{equation}
holds for all $x$ in a cube $E\subseteq \Omega$, where $\Omega$ is a set which can be expressed as the union of rotations and translations of a fixed cube of side $b_{0}$, provided that (a)$E$ has side $b\geq b_{0}$ and (b)$0<d\leq d_{0}$. Here, $\alpha_{n}$ denotes the volume of the unit ball in $R^{n}$, $c$ was defined in (1) and $\Delta_{0}$ was defined in Definition1.3.

The numbers $d_{0}$ and $\omega$ can be expressed specifically as
$$d_{0}=\frac{1}{6C\gamma_{n}(m+1)},\ \omega=\left( \frac{2}{3}\right) ^{\frac{1}{6C\gamma_{n}}}$$
where
$$C=\max\left\{ 2\rho'\sqrt{n}e^{2n\gamma_{n}},\ \frac{2}{3b_{0}}\right\} ,\ \rho'=\frac{\rho}{c}.$$
The numbers $\gamma_{n}$ and $\rho$ were defined in Definition1.1 and 1.3, respectively, and $m=1+\frac{\lambda}{2}$ as in (1).
\end{thm}
{\bf Remark}: (a)The seemingly complicated theorem is in fact not difficult to understand. We expressed each constant specifically for later use. (b)The set $X$ of centers is usually contained in the cube $E$. However it's harmless to extend it to a more general form as a subset of $R^{n}$. (c)$\|f\|_{h}$ is the semi-norm of $f$ mentioned in subsection1.1.

On the right-hand side of the inequality (5), the only things dependent of the shape parameter $c$ are $c^{\frac{\lambda}{2}},\ \omega^{\frac{1}{d}}$ and $\|f\|_{h}$. Also, $d_{0}$ is dependent of $c$. It's not hard to imagine a good choice of $c$ will minimize the error bound (5). However it's nontrivial to express this error bound by an explicit function of $c$, and clarify the relation between $d_{0}$ and $c$. Further treatment of (5) is necessary.
\begin{lem}
Let $\sigma>0$ and $h$ be as in (1). For any $f\in B_{\sigma}$,
\begin{equation}
\|f\|_{h}\leq C_{0}(m,n)\cdot \left(\frac{2}{\pi}\right)^{\frac{1}{4}}\cdot \sigma^{\frac{1+n+\lambda}{4}}\cdot c^{\frac{1-n-\lambda}{4}}\cdot e^{\frac{c\sigma}{2}}\cdot\|f\|_{L^{2}}
\end{equation}
, where
$$C_{0}(m,n)=\frac{(2\pi)^{-n}\sqrt{m!}}{\sqrt{l(\lambda,n)}}.$$
\end{lem}
{\bf Proof.} By Corollary3.3 of \cite{MN2},
\begin{eqnarray*}
\|f\|_{h} & = & \left\{ \sum_{|\alpha|=m}\frac{m!}{\alpha!}\|(D^{\alpha}f)^{\hat{}}\|^{2}_{L^{2}(\rho)}\right\} ^{1/2}\\
          & = & \left\{ \sum_{|\alpha|=m}\frac{m!}{\alpha!}\int |(D^{\alpha}f)^{\hat{}}(\xi)|^{2}\cdot \frac{1}{(2\pi)^{2n}|\xi|^{2m}\hat{h}(\xi)}d\xi \right\} ^{1/2}\\
          & = & (2\pi)^{-n}\left\{ \sum_{|\alpha|=m}\frac{m!}{\alpha!}\int \frac{\xi^{2\alpha}|\hat{f}(\xi)|^{2}}{|\xi|^{2m}l(\lambda,n)|\xi|^{-\lambda-n}\tilde{{\cal K}}_{\frac{n+\lambda}{2}}(c|\xi|)}d\xi\right\} ^{1/2}\\
          & \leq & \frac{(2\pi)^{-n}\sqrt{m!}}{\sqrt{l(\lambda,n)}}\cdot \left\{\int\frac{|\hat{f}(\xi)|^{2}}{|\xi|^{-\lambda-n}\tilde{{\cal K}}_{\frac{n+\lambda}{2}}(c|\xi|)}d\xi\right\}^{1/2}\\
          & = & \frac{(2\pi)^{-n}\sqrt{m!}}{\sqrt{l(\lambda,n)}}\cdot \left\{\int\frac{|\hat{f}(\xi)|^{2}}{|\xi|^{-\lambda-n}(c|\xi|)^{\frac{n+\lambda}{2}}{\cal K}_{\frac{n+\lambda}{2}}(c|\xi|)}d\xi\right\}^{1/2}\\
          & \leq & C_{0}(m,n)c^{\frac{-n-\lambda}{4}}\left\{ \int |\hat{f}(\xi)|^{2}|\xi|^{\frac{n+\lambda}{2}}\cdot \frac{1}{\sqrt{\frac{\pi}{2}}\cdot \frac{e^{-c|\xi|}}{\sqrt{c|\xi|}}}d\xi\right\}^{1/2}\ \mbox{by}\ \cite{We}\ \mbox{where} \\
          &   & C_{0}(m,n):=\frac{(2\pi)^{-n}\sqrt{m!}}{\sqrt{l(\lambda,n)}}\\
          & = & C_{0}(m,n)\left(\frac{2}{\pi}\right)^{\frac{1}{4}}c^{\frac{1-n-\lambda}{4}}\left\{ \int |\hat{f}(\xi)|^{2}|\xi|^{\frac{1+n+\lambda}{2}}e^{c|\xi|}d\xi\right\}^{1/2}. 
\end{eqnarray*}
If $f\in B_{\sigma}$, then 
$$\|f\|_{h}\leq C_{0}(m,n)\left(\frac{2}{\pi}\right)^{\frac{1}{4}}\sigma^{\frac{1+n+\lambda}{4}}c^{\frac{1-n-\lambda}{4}}e^{\frac{c\sigma}{2}}\|f\|_{L^{2}}.$$
\hspace{15cm} $\sharp$
\begin{cor}
Let $\sigma>0$. If $f\in B_{\sigma}$, then (5) can be transformed into
\begin{equation}
|f(x)-s(x)|\leq \sqrt{m!}(2\pi)^{-n}\sqrt{2n\alpha_{n}}\sigma^{\frac{1+n+\lambda}{4}}\sqrt{\Delta_{0}}c^{\frac{1-n+\lambda}{4}}e^{\frac{c\sigma}{2}}\omega^{\frac{1}{d}}\|f\|_{L^{2}}.
\end{equation}
\end{cor}
{\bf Proof.} $B_{\sigma}\subseteq {\cal C}_{h,m}$ is a simple result of Corollary3.3 of \cite{MN2}. Now (7) is just a combination of (5) and (6). \hspace{14.2cm} $\sharp$\\
\\
On the right-hand side of (7), the only things dependent of $c$ are $c^{\frac{1-n+\lambda}{4}},\ e^{\frac{c\sigma}{2}}$ and $\omega^{\frac{1}{d}}$, where $d$ denotes fill distance. It seems that one can abstract a function of $c$ from the long expression of (7) after every thing independent of $c$ is fixed, including $d$.

However as $c$ changes, the upper bound $d_{0}$ of $d$ changes also. As required by Theorem2.1, $d\leq d_{0}$ where 
$$d_{0}=\frac{1}{6C\gamma_{n}(m+1)}.$$
The number $C\rightarrow \infty$, and hence $d_{0}\rightarrow 0$, as $c\rightarrow 0^{+}$. In order to keep $d\leq d_{0}$, the minimal possible choice of $c$ is
\begin{equation}
c_{0}:=12\rho\sqrt{n}e^{2n\gamma_{n}}\gamma_{n}(m+1)d
\end{equation}
where $d$ is fixed. Therefore, we can choose $c$ from the interval $[c_{0}, \infty)$ only.

There is a logical problem about $d,\ d_{0},\ c$ and $b_{0}$. Before $c$ and $b_{0}$ are given, $d_{0}$ is unknown and we do not know whether $d$ is eligible. This question is not difficult to resolve. For any $b_{0}>0$, we require $d<\frac{1}{6C\gamma_{n}(m+1)}$ where we let $C=\frac{2}{3b_{0}}$ temporarily. Then $d<\frac{b_{0}}{4\gamma_{n}(m+1)}$. Note that $2\rho'\sqrt{n}e^{2n\gamma_{n}}=\frac{2}{3b_{0}}$ iff $c=3b_{0}\rho\sqrt{n}e^{2n\gamma_{n}}$. With the same $b_{0}$ and $c\geq 3b_{0}\rho\sqrt{n}e^{2n\gamma_{n}}$, we have
$$C=\max \left\{ 2\rho'\sqrt{n}e^{2n\gamma_{n}},\ \frac{2}{3b_{0}}\right\} =\frac{2}{3b_{0}}$$
and $$d_{0}=\frac{b_{0}}{4\gamma_{n}(m+1)}.$$
Then $d<d_{0}$ is satisfied.

With the same $b_{0}$ and $c<3b_{0}\rho\sqrt{n}e^{2n\gamma_{n}}$, we have 
$$C=\max \left\{ 2\rho'\sqrt{n}e^{2n\gamma_{n}},\ \frac{2}{3b_{0}}\right\} =2\rho' \sqrt{n}e^{2n\gamma_{n}}$$
and 
$$d_{0}=\frac{1}{6C\gamma_{n}(m+1)}=\frac{c}{12\rho\sqrt{n}e^{2n\gamma_{n}}\gamma_{n}(m+1)}.$$
In order to keep $d<d_{0}$, we require 
$$c\geq c_{0}:=12\rho\sqrt{n}e^{2n\gamma_{n}}\gamma_{n}(m+1)d$$
where $d$ was given above satisfying $d<\frac{b_{0}}{4\gamma_{n}(m+1)}$. Therefore, once the cube side $b_{0}$ is given, we first choose $d<\frac{b_{0}}{4\gamma_{n}(m+1)}$. Then put the restriction $c\in [c_{0}, \infty)$ where $c_{0}=12\rho \sqrt{n}e^{2n\gamma_{n}}\gamma_{n}(m+1)d$. The condition $d<d_{0}$ will always be satisfied.

Theoretically $d$, and hence $c_{0}$, can be arbitrarily small. However practically the problem of ill-conditioning when constructing $s(x)$ has to be considered. In this paper we explore theoretically the optimal choice of $c$ for $c_{0}\leq c<\infty$.
\section{Criteria of Choosing $c$}
In the long expression of (7) the only things dependent of $c$ are $c^{\frac{1-n+\lambda}{4}},\ e^{\frac{c\sigma}{2}}$ and $\omega^{\frac{1}{d}}$. Let's define
\begin{equation}
MN(c):=c^{\frac{1-n+\lambda}{4}}\cdot e^{\frac{c\sigma}{2}}\cdot \omega^{\frac{1}{d}}
\end{equation}
and call it an MN function. The optimal choice of $c$ will then be the number minimizing $MN(c)$, when every thing independent of $c$ is fixed.

The value of $\omega$ highly depends on the cube side $b_{0}$. Theoretically, $\omega\rightarrow 0$ and is very influential as $c\rightarrow \infty$ and $b_{0}\rightarrow \infty$ simultaneously. However, a lot of time $b_{0}$ is fixed and cannot approach $\infty$. Therefore we divide our criteria into two classes.
\subsection{$b_{0}$ fixed}
Recall that $C=\max \{ 2\rho'\sqrt{n}e^{2n\gamma_{n}},\ \frac{2}{3b_{0}}\}$ and $\rho'=\frac{\rho}{c}$. The two values $2\rho'\sqrt{n}e^{2n\gamma_{n}}$ and $\frac{2}{3b_{0}}$ are equal if and only if 
\begin{equation}
c=3b_{0}\rho\sqrt{n}e^{2n\gamma_{n}}=:c_{1}
\end{equation}
Then, for $c\in [c_{0}, c_{1}),\ C=2\rho'\sqrt{n}e^{2n\gamma_{n}}$. For $c\in [c_{1}, \infty),\ C=\frac{2}{3b_{0}}$. So we have
\begin{eqnarray*}
\omega^{\frac{1}{d}} & = & \left(\frac{2}{3}\right)^{\frac{1}{6C\gamma_{n}d}}\\
                     & = & \left\{ \begin{array}{ll}
                                     \left(\frac{2}{3}\right)^{\frac{c}{12\rho\sqrt{n}e^{2n\gamma_{n}}\gamma_{n}d}} & \mbox{if $c\in [c_{0},c_{1})$,}\\
                                     \left(\frac{2}{3}\right)^{\frac{b_{0}}{4\gamma_{n}d}} & \mbox{if $c\in [c_{1}, \infty)$}
                                   \end{array} \right. \\
                     & = & \left\{ \begin{array}{ll}
                                     e^{\frac{c\cdot ln\frac{2}{3}}{12\rho\sqrt{n}e^{2n\gamma_{n}}\gamma_{n}d}} & \mbox{\ \ \ if $c\in [c_{0}, c_{1})$,}\\
                                     e^{\frac{b_{0}ln\frac{2}{3}}{4\gamma_{n}d}} & \mbox{\ \ \  if $c\in [c_{1},\infty).$}                
                                   \end{array} \right.   
\end{eqnarray*}
Putting this result into $MN(c)$, we thus have
\begin{equation}
  MN(c)=\left\{ \begin{array}{ll}
                  c^{\frac{1-n+\lambda}{4}}\cdot e^{c\left(\frac{\sigma}{2}+\frac{ln\frac{2}{3}}{12\rho\sqrt{n}e^{2n\gamma_{n}}\gamma_{n}d}\right)} & \mbox{if $c\in [c_{0},c_{1}),$}\\
                  c^{\frac{1-n+\lambda}{4}}\cdot e^{\frac{c\sigma}{2}}\cdot e^{\frac{b_{0}ln\frac{2}{3}}{4\gamma_{n}d}} & \mbox{if $c\in [c_{1},\infty)$}                 \end{array}\right.
\end{equation}
which is a continuous function. Our goal is to find $c$ which minimizes $MN(c)$.

As for the fill distance $d$, we require $d<\frac{b_{0}}{4\gamma_{n}(m+1)}$ once the cube side $b_{0}$ is given, as mentioned in the end of section2.

We have the following cases where $n,\ \lambda,\ m,\ \sigma,\ \gamma_{n},\ \rho$ and $d$ were defined in (1), Definition1.1, 1.2, 1.3 and Theorem 2.1. Moreover we introduce a constant
$$k:=\frac{\sigma}{2}+\frac{ln\frac{2}{3}}{12\rho\sqrt{n}e^{2n\gamma_{n}}\gamma_{n}d}.$$\\
\\
{\bf Case1.} \fbox{$1-n+\lambda>0$} and \fbox{$k\geq 0$}
Let $f\in B_{\sigma}$ be the approximated function, $b_{0}>0$ be the cube side as in Theorem2.1, and $d<\frac{b_{0}}{4\gamma_{n}(m+1)}$ be the fill distance. If $1-n+\lambda>0$ and $k\geq 0$, the optimal choice of $c$ for $c\in [c_{0}, \infty)$ is \fbox{$c=12\rho\sqrt{n}e^{2n\gamma_{n}}\gamma_{n}(m+1)d$}.\\
\\
{\bf Reason}: In this case $MN(c)$ is an increasing function in its domain $[c_{0},\infty)$. We therefore choose $c=c_{0}$, defined in (8), which minimizes $MN(c)$ in (11) and the error bound in (7).\\
\\
{\bf Case2.} \fbox{$1-n+\lambda>0$} and \fbox{$k<0$} Let $f\in B_{\sigma}$ be the approximated function, $b_{0}>0$ be the cube side as in Theorem2.1, and $d<\frac{b_{0}}{4\gamma_{n}(m+1)}$ be the fill distance. If $1-n+\lambda>0$ and $k<0$, the optimal choice of $c$ for $c\in [c_{0},\infty)$ is the value minimizing $MN(c)$ over the interval $[c_{0},c_{1}]$ where $c_{1}$ was defined in (10).\\
\\
{\bf Reason.} In this case $MN(c)$ is increasing on $[c_{1},\infty)$. Therefore the minimum value of $MN(c)$ in $[c_{0},\infty)$ must happen in $[c_{0},c_{1}]$.\\
\\
{\bf Remark}:(a)This case rarely happens because $k$ is usually positive. (b)In fact, the optimal $c$ in Case2 can be obtained exactly. Let $g(c):=c^{\frac{1-n+\lambda}{4}}e^{ck}$ which is just $MN(c)$ on $[c_{0},c_{1}]$. Then $g(c)$ is increasing on $(0,\frac{1-n+\lambda}{-4k}]$ and decreasing on $[\frac{1-n+\lambda}{-4k},\infty)$. So the optimal $c$ is (i)$c_{0}$ if $c_{1}\leq \frac{1-n+\lambda}{-4k}$, (ii)$c_{1}$ if $\frac{1-n+\lambda}{-4k}\leq c_{0}$, and (iii)$c_{1}$ if $g(c_{1})\leq g(c_{0})$, and $c_{0}$ if $g(c_{0})\leq g(c_{1})$, for $c_{0}<\frac{1-n+\lambda}{-4k}<c_{1}$.\\
\\
{\bf Case3.} \fbox{$1-n+\lambda<0$} and \fbox{$k=0$} \   Let $f\in B_{\sigma}$ be the approximated function, $b_{0}>0$ be the cube side as in Theorem2.1, and $d<\frac{b_{0}}{4\gamma_{n}(m+1)}$ be the fill distance. If $1-n+\lambda<0$ and $k=0$, the optimal choice of $c$ for $c\in [c_{0},\infty)$ is the value minimizing $c^{\frac{1-n+\lambda}{4}}e^{\frac{c\sigma}{2}}$ on the interval $[c_{1}, \infty)$.\\
\\
{\bf Reason}: In this case $MN(c)$ is decreasing on $[c_{0},c_{1})$. Therefore the minimum value of $MN(c)$ happens on the interval $[c_{1},\infty)$ where the essential part of $MN(c)$ is $c^{\frac{1-n+\lambda}{4}}e^{\frac{c\sigma}{2}}$.\\
\\
{\bf Remark}:(a)Although $k$ rarely equals zero, we can make it zero by choosing $\sigma$ in an appropriate way. It will make things easier for the optimal choice of $c$. (b)If $1-n+\lambda<0$, it can be shown that $g(c):=c^{\frac{1-n+\lambda}{4}}e^{\frac{c\sigma}{2}}$ is decreasing on $(0, \frac{-1+n-\lambda}{2\sigma}]$ and increasing on $[\frac{-1+n-\lambda}{2\sigma},\infty)$, with $g'(\frac{-1+n-\lambda}{2\sigma})=0$. Therefore the optimal choice of $c$ for $c\in [c_{0},\infty)$ in Case3 is in fact \fbox{$\max\left\{ \frac{-1+n-\lambda}{2\sigma},c_{1}\right\}$}.\\
\\
{\bf Case4}. \fbox{$1-n+\lambda<0$} and \fbox{$k>0$}\ Let $f\in B_{\sigma}$ be the approximated function, $b_{0}>0$ be the cube side as in Theorem2.1, and $d<\frac{b_{0}}{4\gamma_{n}(m+1)}$ be the fill distance. Let $g_{1}(c):=MN(c)|_{c\in[c_{0},c_{1})}$ and $g_{2}(c):=MN(c)|_{c\in [c_{1},\infty)}$. If $1-n+\lambda<0$ and $k>0$, then
\begin{list}
  {(\roman{milk})}{\usecounter{milk} \setlength{ \rightmargin}{\leftmargin}}
  \item $g_{1}'(\frac{1-n+\lambda}{-4k})=0,\ g_{1}'(c)<0$ for $c\in (0,\frac{1-n+\lambda}{-4k})$, and $g_{1}'(c)>0$ for $c\in (\frac{1-n+\lambda}{-4k},\infty)$, and
  \item $g_{2}'(\frac{-1+n-\lambda}{2\sigma})=0,\ g_{2}'(c)<0$ for $c\in(0,\frac{-1+n-\lambda}{2\sigma})$, and $g_{2}'(c)>0$ for $c\in (\frac{-1+n-\lambda}{2\sigma},\infty)$. 
\end{list}
In this case the minimum value of $g_{1}(c)$ on $[c_{0},c_{1})$ happens at (a)$c=c_{0}$ if $\frac{1-n+\lambda}{-4k}\leq c_{0}$, (b)$c=\frac{1-n+\lambda}{-4k}$ if $c_{0}<\frac{1-n+\lambda}{-4k}<c_{1}$, and (c)$c=c_{1}$ if $c_{1}\leq \frac{1-n+\lambda}{-4k}$. Also, the minimum value of $g_{2}(c)$ on $[c_{1},\infty)$ happens at (a)$c=c_{1}$ if $\frac{-1+n-\lambda}{2\sigma}<c_{1}$, and (b)$c=\frac{-1+n-\lambda}{2\sigma}$ if $c_{1}\leq \frac{-1+n-\lambda}{2\sigma}$.

Let $c^{*}\in [c_{0},c_{1})$ minimize $g_{1}(c)$ and $c^{**}\in [c_{1},\infty)$ minimize $g_{2}(c)$. Then the optimal choice of $c\in [c_{0},\infty)$ is (a)$c^{*}$ if $g_{1}(c^{*})\leq g_{2}(c^{**})$, and (b)$c^{**}$ if $g_{2}(c^{**})\leq g_{1}(c^{*})$.\\
\\
{\bf Reason}: By direct differentiation, we get
$$g_{1}'(c)=c^{\frac{1-n+\lambda}{4}-1}e^{ck}\left(\frac{1-n+\lambda}{4}+ck\right)$$
and 
$$g_{2}'(c)=c^{\frac{1-n+\lambda}{4}-1}e^{\frac{c\sigma}{2}}e^{\frac{b_{0}ln\frac{2}{3}}{4\gamma_{n}d}}\left(\frac{1-n+\lambda}{4}+\frac{c\sigma}{2}\right).$$
The seemingly complicated criterion is then just a simple result of the two derivatives. \hspace{2cm} $\sharp$\\
\\
{\bf Case5}. \fbox{$1-n+\lambda<0$} and \fbox{$k<0$} Let $f\in B_{\sigma}$ be the approximated function, $b_{0}>0$ be the cube side as in Theorem2.1, and $d<\frac{b_{0}}{4\gamma_{n}(m+1)}$ be the fill distance. If $1-n+\lambda<0$ and $k<0$, the optimal choice of $c$ for $c\in [c_{0},\infty)$ is (a)$c=c_{1}$ if $\frac{-1+n-\lambda}{2\sigma}<c_{1}$, and (b)$c=\frac{-1+n-\lambda}{2\sigma}$ if $c_{1}\leq \frac{-1+n-\lambda}{2\sigma}$.\\
\\
{\bf Reason}: In this case $MN(c)$ is decreasing on $[c_{0},c_{1})$ and its minimum value happens at a number in $[c_{1},\infty)$ which minimizes $g_{2}$ of Case4. Our criterion thus follows immediately from Case4. \hspace{1.1cm} $\sharp$ 
\subsection{$b_{0}$ not fixed}
In Theorem2.1 one can easily find that if $b_{0}$ is not fixed and can be made arbitrarily large, $C$ will become arbitrarily small by letting $c\rightarrow \infty$. This will make $\omega \rightarrow 0$ and $d_{0}\rightarrow \infty$, a very beneficial situation. In fact some domain $\Omega$ does allow $b_{0}\rightarrow \infty$. For example $\Omega=R^{n}$ or
$$\Omega=\{ (x_{1},\cdots,x_{n})|\ x_{i}>0\ \mbox{for}\ i=1,2,\cdots,n\}.$$
This kind of domain is called {\bf dilation-invariant} by Madych in \cite{MN3}. For this kind of domain we have a different set of criteria of choosing $c$.

By increasing the cube side $b_{0}$, one can always keep $C=2\rho'\sqrt{n}e^{2n\gamma_{n}}$ in Theorem2.1. Thus 
$$MN(c)=c^{\frac{1-n+\lambda}{4}}e^{c\left(\frac{\sigma}{2}+\frac{ln\frac{2}{3}}{12\rho\sqrt{n}e^{2n\gamma_{n}}\gamma_{n}d}\right)}$$
for all $c\in [c_{0}, \infty)$.

Note that we never decrease $b_{0}$ because it will make the error bound (7) worse.

As for the choice of $d$, there is no restriction. Once $d>0$ is given, $c_{0}$ will be fixed. For any $c\in [c_{0},\infty)$, the condition $d<d_{0}$ in Theorem2.1 will be satisfied. However the smaller $d$ is, the larger $[c_{0},\infty)$ is, making the criteria more meaningful.

Once $c$ is chosen, we let $b_{0}=\frac{c}{3\rho\sqrt{n}e^{2n\gamma_{n}}}$ such that $2\rho'\sqrt{n}e^{2n\gamma_{n}}=\frac{2}{3b_{0}}$ in the definition of $C$. Further increasing $b_{0}$ is not to be expected because more data points will be involved to keep the fill distance $d$ fixed.

Note that by letting $k:=\frac{\sigma}{2}+\frac{ln\frac{2}{3}}{12\rho\sqrt{n}e^{2n\gamma_{n}}d}$, we have a more simple expression for $MN(c)$ as $$MN(c)=c^{\frac{1-n+\lambda}{4}}e^{ck}.$$

The way of choosing $c$ is then divided in the following cases.\\
\\
{\bf Case1.} \fbox{$1-n+\lambda>0$} and \fbox{$k\geq 0$} Let $f\in B_{\sigma}$ be the approximated function, and the domain $\Omega$ in Theorem2.1 be dilation-invariant. For any fill distance $d>0$, if $1-n+\lambda>0$ and $k\geq 0$, the optimal choice of $c$ for $c\in [c_{0},\infty)$ is $c=c_{0}$.\\
\\
{\bf Reason}: In this case $MN(c)$ is increasing on $[c_{0},\infty)$. \hspace{6.8cm} $\sharp$\\
\\
{\bf Case2}. \fbox{$1-n+\lambda>0$} and \fbox{$k<0$} Let $f\in B_{\sigma}$ be the approximated function, and the domain $\Omega$ in Theorem2.1 be dilation-invariant. For any fill distance $d>0$, if $1-n+\lambda>0$ and $k<0$, the larger $c$ is, the better it is.\\
\\
{\bf Reason}: In this situation,
\begin{eqnarray*}
  MN'(c) & = & c^{\frac{1-n+\lambda}{4}-1}e^{ck}\left(\frac{1-n+\lambda}{4}+kc\right)\\
         & = & 0
\end{eqnarray*}
iff
$$c=\frac{1-n+\lambda}{-4k}.$$
The function $MN(c)$ is increasing and decreasing, respectively, depending on $c<\frac{1-n+\lambda}{-4k}$ or $c>\frac{1-n+\lambda}{-4k}$. Note that $MN(c_{0})$ is a finite value and $MN(c)\rightarrow 0$ as $c\rightarrow \infty$. Our criterion thus follows.$\sharp$\\
\\
{\bf Remark}:(a)Usually $k\geq 0$. Hence Case2 rarely happens. (b)If $c_{0}<\frac{1-n+\lambda}{-4k}$, $MN(c_{0})$ may be small enough to be accepted. If it happens, $c=c_{0}$ is also a good choice. \\
\\
{\bf Case3}. \fbox{$1-n+\lambda<0$} and \fbox{$k>0$} Let $f\in B_{\sigma}$ be the approximated function, and the domain $\Omega$ in Theorem2.1 be dilation-invariant. For any fill distance $d>0$, if $1-n+\lambda<0$ and $k>0$, the optimal choice of $c$ for $c\in [c_{0},\infty)$ is (a)$c=\frac{1-n+\lambda}{-4k}$ if $c_{0}\leq \frac{1-n+\lambda}{-4k}$ and (b)$c=c_{0}$ if $c_{0}>\frac{1-n+\lambda}{-4k}$. \\
\\
{\bf Reason}: By the structure of $MN'(c)$, one sees easily $MN(c)$ is decreasing on $(0,\frac{1-n+\lambda}{-4k}]$ and increasing on $[\frac{1-n+\lambda}{-4k},\infty)$. Our criterion thus follows immediately. \hspace{5cm} $\sharp$\\
\\
{\bf Case4}: \fbox{$1-n+\lambda<0$} and \fbox{$k\leq 0$} Let $f\in B_{\sigma}$ be the approximated function, and the domain $\Omega$ in Theorem2.1 be dilation-invariant. For any fill distance $d>0$, if $1-n+\lambda<0$ and $k\leq 0$, the larger $c$ is, the better it is.\\
\\
{\bf Reason}: In this case $MN(c)$ is decreasing on $(0,\infty)$. Moreover, $MN(c)\rightarrow 0$ as $c\rightarrow \infty$. \hspace{1.6cm} $\sharp$\\
\\
{\bf Remark}:(a)Case1-4 are all based on Theorem2.1. Although $c$ and $b_{0}$ appear before $d$ in Theorem2.1, in the four cases $d$ is given before $c$ and $b_{0}$ are chosen. Once $d$ is fixed, we chose $c$ according to the four criteria. After $c$ is chosen, we let $b_{0}=\frac{c}{3\rho\sqrt{n}e^{2n\gamma_{n}}}$.(b)In order to keep the fill distance $d$ fixed, more data points have to be added when the cube side $b_{0}$ increases. (c)Although in this paper we only explore $c$ over the interval $[c_{0},\infty)$, the curve of $MN(c)$ can be used as a conjecture for $c\in (0,c_{0}]$. It probably holds even if we cannot prove it.


\end{document}